\documentclass[10pt]{article}
\usepackage{amsfonts}

\newcommand{\lar}{\longrightarrow}
\newcommand{\kA}{\mathcal A}
\newcommand{\kE}{\mathcal E}
\newcommand{\kF}{\mathcal F}
\newcommand{\kQ}{\mathcal Q}
\newcommand{\kK}{\mathcal K}
\newcommand{\kG}{\mathcal G}
\newcommand{\kH}{\mathcal H}
\newcommand{\kL}{\mathcal L}
\newcommand{\kB}{\mathcal B}
\newcommand{\kM}{\mathcal M}
\newcommand{\kO}{\mathcal O}
\newcommand{\kP}{\mathcal P}
\newcommand{\kI}{\mathcal I}
\newcommand{\kC}{\mathcal C}
\newcommand{\qed}{{\begin{flushright}$\diamond$\end{flushright}}}
\renewcommand{\P}{\mathbb P}
\newcommand{\A}{\mathbb A}

\newtheorem{theorem}{Theorem}
\newtheorem{corollary}{Corollary}

\newtheorem{definition}{Definition}

\begin{document}
\title{Geometrical Properties of Sections of Buchsbaum-Rim Sheaves\\
       or\\
       How to Construct Gorenstein Schemes of Higher Codimension with
       {\sl Singular}}
\author{Igor Burban \and Hans Georg Freiermuth\thanks{
        Both authors would like to thank the Universita di Catania
        and the Universit\"at Kaisers\-lautern for the generous support.}}
\date{}
\maketitle

\begin{center} {\bf 1.~Introduction} \end{center}
In this {\it survey article} we want to discuss one of the ways to construct
arithmetically Gorenstein subvarieties of projective space.
It is well-known that in codimension $2$ an arithmetically Gorenstein
subvariety is always a complete intersection.  For codimension $\geq
4$, however, the construction of Gorenstein subschemes, apart
from complete intersections, is quite a complicated problem since
no structure theorem is known as in codimension $3$ (\cite{Eisenbud}).
On the other hand, it is necessary to develop a technique for
constructing such schemes in view of Gorenstein liaison.\\
Consider the reflexive kernel sheaves $\kB_\phi$ of sufficiently
general, generically surjective morphisms $\phi$ between decomposable bundles 
over $\P^n$, so called Buchsbaum-Rim
sheaves. The desired Gorenstein schemes appear quite unexpectedly as the
top-dimensional part of the zero-locus of a regular section
$s\in H^0(\P^n,\kB_\phi)$ (\,\cite{Codim3}, \cite{Multsect}).
As an application one gets information about the geometrical
properties of sections of certain non-split rank $n$ vector
bundles on $\P^n$.\\
We show a way how to implement the construction method in the
computer algebra system {\it Singular} and produce some examples of
Gorenstein curves and threefolds in $\P^6$.\\
Finally, we investigate a class of sheaves $\kB_\phi$ where the
degeneracy locus of $\phi$ does not have the expected codimension.

\begin{center} {\bf 2.~Gorenstein points in $\P^3$} \end{center}
Let $R=k[z_0,z_1,z_2,z_3]$, $k=\bar{k}$ and $A$ a homogeneous $t \times
(t+3)$-matrix over $R$ such that ideal of all $t\times t$-minors of
$A$ is $(z_0,z_1,z_2,z_3)$-primary. $A$ is given as a graded degree
zero homomorphism between two free $R$-modules $F$ and $G$ of rank
$t+3$ and $t$ respectively:
$$
0 \lar Q \lar F \stackrel{A}\lar G \lar coker(A) \lar 0.
$$
Because of our conditions on the minors, $coker(A)$ has finite length.
Sheafifying the sequence, one gets:
$$
0 \lar \kQ \lar \kF \lar \kG \lar 0.
$$
Then $\kQ$ is a vector bundle of rank $3$. We can assume that $\kQ$
has global sections - if necessary twist it. Choose a regular global
section $s$, i.e.~one with a zero-dimensional degeneracy locus.
\begin{theorem}[\cite{Codim3}]\label{points}
Let $X=Z(s)$ be the zero locus of $s$. Then $X$ is arithmetically
Gorenstein and its saturated ideal $I_X$ has a free resolution
$$
0 \lar R(-c_1) \lar F(-c_1) \oplus G^{\ast} \lar G(-c_1)\oplus F^{\ast}
\lar I_X \lar 0,
$$
where $c_1=c_1(\kQ)$ denotes the first Chern class of $\kQ$.
\end{theorem}
{\bf Idea of the proof:} The section $s$ determines an exact sequence
$$
0\lar \kO_{\P^3} \stackrel{s}\lar \kQ \lar \kC \lar 0.
$$
One can easily show that the cokernel $\kC$ is reflexive of rank $2$.
Hence $\kC^{\ast}=\kC(-c_1)$ where $c_1(\kQ)=c_1(\kC)$ is the first Chern
class. Dualizing the above sequence, we get
$$
0\lar \kC(-c_1) \lar \kQ^{\ast} \stackrel{s^t}\lar \kO_{\P^3}\lar 
\kO_{X} \lar 0,
$$
since ${\mathcal Ext}^{1}(\kC, \kO_{\P^3})=\kO_{X}$ and $\kQ$ is locally free.
Using $H^1_\ast(\P^3,\kC)=H^1_\ast(\P^3,\kQ)$ one obtains a diagram
\[\begin{array}{ccccccc}
   0 & & 0 & & & & 0 \\
  \uparrow & & \uparrow & & & & \uparrow \\
0  \lar  H^0_\ast(\kC(-c_1)) & \lar  & H^0_\ast(\kQ^\ast) & \lar & I_X &
    \lar & H^1_\ast(\kC(-c_1)) \lar 0\\
   \uparrow & & \uparrow & & & & \uparrow  \\
   E_1(-c_1) & & F^\ast & & & & G(-c_1)  \\
   \uparrow & & \uparrow & & & & \uparrow  \\
   R(-c_1)\oplus E_2(-c_1) & & G^\ast & & & & F(-c_1) \\
   \uparrow & & \uparrow & & & & \uparrow  \\
   E_3(-c_1) & & 0 & & & & E_1(-c_1)\\
   \uparrow & & & & & & \uparrow\\
    0       & & & & & & E_2(-c_1)\\
            & & & & & & \uparrow\\
            & & & & & & E_3(-c_1)\\
            & & & & & & \uparrow\\
                  & & & & & & 0
\end{array}\]
where $E_i=\wedge^{t+i} F\otimes S\,^{i-1}(G)^\ast\otimes
\wedge^t\,G^\ast$ are the $R$-modules from the acyclic (!)
Buchsbaum-Rim complex associated to $F\stackrel{A}{\lar} G$.
Now the application of a mapping cone, the horseshoe lemma and
$\mbox{hd}_R I_X=2$ (X is arithmetically Cohen-Macaulay) yield an 
exact sequence
\begin{eqnarray*}
 0\lar R(-c_1)\oplus E_1(-c_1) \lar F(-c_1)\oplus G^\ast \oplus
E_1(-c_1)\lar G(-c_1)\oplus F^\ast \lar & & \\
\cdots\lar I_X \lar 0 & &
\end{eqnarray*}
Using this sequence, one proves immediately that the $h$-vector
of $R/I_X$ is symmetric. Furthermore, one shows that $X$ has the
generalized Cayley-Bacharach property with respect to
$|\,K_{\P_3} \otimes \deg Q\,|$. But this implies that
$X$ is arithmetically Gorenstein and therefore
we can delete $E_1(-c_1)$ in the above resolution.
\qed
\begin{corollary}[\cite{Codim3}]
Let $s$ be a regular global section of a rank three vector bundle $\kE$ on
${\P^3}$ with $H^{2}_{*}({\P^3},\kE)=0$ but non-vanishing
first cohomology. Then the zero scheme $X=Z(s)$ is arithmetically Gorenstein.
\end{corollary}
{\bf Proof:}  The module $H^{1}_{*}({\P^3},\kE)$ has finite
length by the Enriques-Zariski-Severi vanishing lemma and thus the
Auslander-Buchsbaum theorem implies the existence of a minimal free
resolution
$$
0 \lar F_4 \lar F_3 \lar F_2 \lar F \lar G \lar H^{1}_{*}(
\P^3,\kE)
\lar 0.
$$
After sheafifying, we get
$$
0 \lar \kF_4 \lar \kF_3 \lar \kF_2 \lar \kF \stackrel{\phi}\lar \kG
\lar 0.
$$
From the sequence
$$
0 \lar ker(\phi) \lar \kF \lar \kG \lar 0
$$
and the associated long exact sequence in cohomology
we deduce that $ker(\phi)$ is a vector bundle with
$H^{1}_{*}(ker(\phi))\cong H^{1}_{*}(\kE)$,
$H^{2}_{*}(ker(\phi))=0$. Thus, $\kE=ker(\phi)\oplus \kL$, where $\kL$
is a decomposable bundle. But $\kE$ does not split by Horrock's theorem
and therefore $\kE=ker(\phi)$. Now apply Theorem \ref{points}.
\qed
{\bf Example 1:} Consider the $2\times 5$-matrix
$$\phi=\left(
\begin{array}{ccccc}
z_0^3+z_1^3 & z_0^3+z_2^3 & z_0^3+z_3^3 & z_1^3+z_2^3 & z_1^3+z_3^3\\
z_2^3+z_3^3 & z_0^3       & z_1^3       & z_2^3       & z_3^3
\end{array}
 \right)$$
as map in the exact sequence
$$0\lar\kQ\lar 5\,\kO_{\P^3}(6)\stackrel{\phi}{\lar}2\,\kO_{\P^3}(9)\lar 0.$$
Here $c_1(\kQ)=12$. The ideal $I$ of the vanishing locus of
a randomly chosen regular section $s\in H^0(\P^3,\kQ)$ is for example
$$\begin{array}{rl}
 J= &(\, 103\,z_1^3\,z_2^3+66\,z_2^6-763\,z_0^3\,z_3^3+660\,z_1^3\,z_3^3+
  829\,z_2^3\,z_3^3-713\,z_3^6, \\
  & 618\,z_0^3\,z_2^3+90\,z_2^6-5741\,z_0^3z_3^3+5123\,z_1^3\,z_3^3+7067\,
  z_2^3z_3^3-5935\,z_3^6, \\
  &
  103\,z_1^6-168\,z_2^6+959\,z_0^3z_3^3-959\,z_1^3\,z_3^3-1230\,z_2^3\,z_3^3+
  1019\,z_3^6,  \\
  &
  103\,z_0^3\,z_1^3-174\,z_2^6+710\,z_0^3\,z_3^3-710\,z_1^3\,z_3^3-987\,
  z_2^3\,z_3^3+934\,z_3^6, \\
  &  4326\,z_0^6-8814\,z_2^6+12847\,z_0^3\,z_3^3-9139\,z_1^3\,z_3^3-
   13009\,z_2^3\,z_3^3+18923\,z_3^6\, ). 
\end{array}$$
The non-saturated zero-scheme $Z(I)$ of degree $54$ has $h$-vector
$$h_I=(1,3,6,10,15,21,23,21,15,7,-3,-15,-20,-18,-9,-3)$$
and is consequently not Gorenstein.
Computing the saturation $J$ of $I$, we get the Gorenstein ideal
$$(\,13\,z_1^3+18\,z_2^3+31\,z_3^3,\,
   546\,z_0^3+534\,z_2^3+2471\,z_3^3,\,
   5976\,z_2^6-32430\,z_2^3z_3^3-90965\,z_3^6\,)$$
with $\dim Z(J)=0$ and $\deg Z(J)=54$ -- as expected in this
case (cf.~Remark 2).
The symmetric $h$-vector is $h_J=(1,3,6,8,9,9,8,6,3,1)$ and a
minimal free resolution is:
$$0\lar R(-12) \lar R(-6)\oplus 2\,R(-9) \lar 2\,R(-3)\oplus R(-6) \lar I
\lar 0.$$
Thus, we observe that the free resolution in Theorem 1 is not
necessarily minimal.\\[0.2cm]
{\bf Remark 1:} We want to mention that not {\sl all} Gorenstein zero
schemes in $\P^3$ can be obtained in this way. Indeed, the
Buchsbaum-Eisenbud structure theorem \cite{Eisenbud} tells us that every
codimension $3$ arithmetically Gorenstein subscheme in $\P^n$ is precisely the
zero set of the ideal $I$ generated by the $2N\times 2N$
Pfaffians ($=$ roots of the $2N\times 2N$ principal minors) of some
$(2N+1)\times (2N+1)$
skew-symmetric matrix over $R=k[z_0,z_1,\ldots,z_n]$.\\[0.2cm]
Let $R=k[z_0,z_1,z_2,z_3]$ and
$$
A =
\left(
\begin{array}{ccccc}
0 & -z_2 & 0 & z_3 & 0 \\
z_2 & 0 & -z_1 &0 & -z_3 \\
0 & z_1 & 0 &  0 & -z_2 \\
-z_3 & 0 & 0 & 0 & -z_1 \\
0 & z_3 & z_2  & z_1 & 0 \\
\end{array}
\right).
$$
A {\sl minimal} generating system of the ideal $I$ is in this case for
example
$$(z_1^2, z_2^2,z_{1}z_{2}-z_3^2, z_{1}z_{3}, z_{2}z_{3}).$$
It defines a non-reduced Gorenstein point $X$ of degree $5$
and cannot be obtained via a regular section of a rank $3$
Buchsbaum-Rim sheaf. The main reason is that $X$
cannot be generated scheme-theoretically from $4$ equations and
is therefore not obtained from a homomorphism
$$ \phi\,:\,\bigoplus_{i=1}^4 R(a_i)\longrightarrow R(b).$$
The $t\times (t+3)$ matrices with $t>1$ are also out of the 
question for degree reasons as we will see below.\\[0.2cm]
{\bf Remark 2:} Suppose that
$\kF=\bigoplus_{i=1}^{t+3}\kO_{\P_3}(a_i)$, $\kG=\bigoplus_{j=1}^{t}
\kO_{\P_3}(b_j)$ and let $\omega$ denote the hyperplane class.
We want to compute the degree of $X$ in terms of the integers
$a_i$ and $b_j$.
Clearly, the fundamental class $[X]$ equals $c_3(\kQ)$.
The Chern polynomial of the bundle $\kQ$ is
$$c_\omega(\kQ)=c_\omega(\kF)c_\omega(\kG)^{-1}=
\frac{\prod_{i=1}^{t+3} (1+a_i\,\omega)}{\prod_{j=1}^{t}
  (1+b_j\,\omega)}=\sum_{i=0}^3
s_i(a)\,\omega^i\,\left(\sum_{j=0}^{\min(t,3)}s_j(b)\,
   \omega^j\right)^{-1},$$
where $s_i(a)=s_i(a_1,\ldots,a_{t+3})$ and
$s_j(b)=s_j(b_1,\ldots,b_t)$ respectively denote the
elementary symmetric functions.
A straightforward computation gives:
$$\deg(X)=\left\{\scriptsize \begin{array}{l@{,\;}l}
    s_3(a)-b_1\,s_2(a)+b_1^2\,s_1(a)-b_1^3 & t=1\\
    s_3(a)-(b_1+b_2)\,s_2(a)+(b_1^2+b_1b_2+b_2^2)\,s_1(a)-
    (b_1^3+b_1^2b_2+b_1b_2^2+b_2^3) & t=2\\
    s_3(a)-s_3(b)-s_1(a)\,s_2(b)-s_2(a)\,s_1(b)+2s_1(b)\,s_2(b)-
    s_1(b)^2\,s_1(a)-s_1(b)^3 & t\geq 3
    \end{array}\right.$$
Therefore one can use the described method to construct
Gorenstein codimension three schemes of a {\sl given degree}.
\begin{center}
{\bf 3.~The General Case}
\end{center}
Migliore, Nagel and Peterson extended the above result in
\cite{Multsect} to higher, odd codimension.
This provides a construction technique for Gorenstein subschemes
of $\P^n$ with prescribed degree in cases where no structure theorem
is known, as for example in codim $5,7,9,\ldots$\\
One can also arrange that the new subschemes contain a given
equidimensional subscheme of the same codimension. This is very useful
from the viewpoint of Gorenstein liaison.\\[0.2cm]
The setup in the general
situation is the following:\\[0.2cm]
Let $Z=\mbox{Proj}(R)$, where $R$ is a graded Gorenstein $k$-algebra of
$\dim(R)=n+1$. Let $\phi\,:\,\kF \lar \kG$ be a morphism of
vector bundles of rank $f$ and $g$ respectively, $f>g$ such that
\begin{enumerate}
\item the degeneracy locus of $\phi$ has codimension $f-g+1$.
\item $F:=H^{0}_{*}(Z,\kF)$  and $G:=H^{0}_{*}(Z,\kG)$ are free
      $R$-modules.
\end{enumerate}
Consider the exact sequence
$$
0 \lar \kB_{\phi} \lar \kF \stackrel{\phi}\lar \kG \lar \kM_{\phi} \lar
0.
$$
$\kB_{\phi}$ is called a Buchsbaum-Rim sheaf. It is reflexive (as a
second syzygy module) of rank $r=f-g$. Let $\kP$ be a decomposable 
vector bundle of rank $q$, $1\le q < r$, and $\psi : \kP \lar
\kB_{\phi}$ be a ``generalized section''.\\
It induces a map $\wedge^{q}\psi^\ast \,:\,
\wedge^{q}\kB^{\ast}_{\phi}\lar \wedge^{q}\kP^{\ast} $ where
$\wedge^{q}\kP^{\ast}$ is a line bundle.
Let $\delta_{q}:=\wedge^{q}\psi^{\ast}$. Then we get a degeneracy locus $S$
$$
\wedge^{q}\kB^{\ast}_{\phi}\otimes (\wedge^{q}P^\ast)^{-1}
\stackrel{\delta_{q}\otimes id}\lar \kO_{Z} \lar \kO_{S} \lar 0.
$$
We furthermore assume that the generalized section $\psi$ is regular,
i.e.~$S$ has the expected codimension $r-q+1$.
Let $I(\psi)$ denote the saturated ideal of $S$, $X$ the top-dimensional part
of $S$, $J=J(\psi)$ its saturated ideal and $P:=H^0_{\ast}(\P^n,\kP)$.\\[0.2cm]
{\bf The aim of the story} : Study the geometrical properties of $S$
and $X$.\\[0.2cm]
We cite parts of the main result obtained in \cite{Multsect}.
\begin{theorem} Using the notation from above, one has the following:
\begin{itemize}
\item If $r+q$ is odd then $X=S$ is arithmetically Cohen-Macaulay iff
$q=1$. In this case, $X$ has Cohen-Macaulay type $\leq 1+
{\frac{r}{2}+g-1 \choose g-1}$.
\item If $r+q$ is even then $X$ is arithmetically Cohen-Macaulay iff
$q=1,2$.
 In case $q=1$, $X$ is arithmetically Gorenstein. If moreover $r<n$ then
components of $S$ has either codimension
$r-q+1$ or codimension $r+1$.
\item Moreover, there is a resolution
$$
0\lar A_r \oplus C_r \lar \dots \lar A_1\oplus C_1 \lar I_X \otimes
\wedge^{q}P^{\ast} \lar 0,
$$
where
\begin{eqnarray*}
A_{k} & = &
\bigoplus\limits_{\scriptsize
\begin{array}{c}
i+2j=k+q-1  \\ q\le i+j \le \frac{r+q-1}{2}
\end{array}}
 \wedge^{i}F^{\ast}\otimes S^{j}(G)^{\ast}\otimes S^{i+j-q}(P),\\
C_{k} & = &
\bigoplus\limits_{\scriptsize
\begin{array}{c}
i+2j=r+1-q-k  \\ i+j \le \frac{r-q}{2}
\end{array}}
\wedge^{i}F\otimes S^{j}(G) \otimes S^{r-q-i-j}(P)\otimes
\wedge^{f}F^{\ast}\otimes
\wedge^{g}G.
\end{eqnarray*}
\end{itemize}
\end{theorem}
{\bf Important techniques used in the proof:}
Local cohomology; the Eagon-Northcott complex, its dual complex and
their homology; the Buchsbaum-Rim complex.
\qed
{\bf Remark 3:} In order to clear the fog in the jungle of wedges and
symmetric powers in Theorem $2$, we explicitly write down the
resolution of $I_X$ for the case of a regular section in a
rk-$5$ Buchsbaum-Rim sheaf $\kB_\phi$ on $\P^6$
(i.e.~$q=1,r=5$, $R=k[z_0,\ldots,z_6]$ and $c_1:=c_1(\kB_\phi)$):
\begin{eqnarray*}
0\lar R(-c_1)\lar S^2(G)^\ast \oplus F(-c_1) \lar F^\ast\otimes G^\ast
\oplus G(-c_1) \oplus (\wedge^2 F)(-c_1)\lar & & \\
\lar (\wedge^2 F)^{\ast}\oplus G^{\ast} \oplus F\otimes G(-c_1)\lar F^\ast
\oplus S^2(G)(-c_1)\lar I_X\lar 0 & &
\end{eqnarray*}
\begin{corollary}
In the assumptions of the previous theorem, let
$R=k[z_0,z_1,\dots,z_n]$, $q=1$ and $r$ be an odd integer.
Choose a homogeneous $t\times (t+r)$-matrix $A$ over $R$, defining a
morphism
between free modules $F \stackrel{A}\lar G$. Suppose that the ideal
of all $t\times t$-minors of $A$ has the expected codimension $r+1$.
After sheafifying, we get an exact sequence
$$\begin{array}{ccccc}
0\lar \kB_A \lar & \kF & \stackrel{A}\lar & \kG & \lar coker({A}) \lar
0.\\
 & || & & || & \\
 & \bigoplus_{i=1}^{t+r}\kO_{\P^n}(a_i)  & \stackrel{A}\lar &
\bigoplus_{j=1}^t\kO_{\P^n}(b_j) &
\end{array}
$$
Let $s$ be a regular section of $\kB_A$, $S:=Z(s)$ and $X$ be the
top-dimensional part of $S$. Then $X$ is arithmetically Gorenstein of
codimension $r$.
\end{corollary}
Using the Buchsbaum-Rim resolution of $\kB_A$, one 
can compute $\deg(X)=c_r(\kB_A)$ in terms of the twists $a_i$ and
$b_j$ of $\kF$ and $\kG$ respectively in a similar fashion as for the
Gorenstein points in $\P^3$.\\[0.2cm]
{\bf Remark 4:} It can be shown that the Buchsbaum-Rim sheaves on $\P^n$
are exactly the reflexive Eilenberg-MacLane sheaves $\kB$ of rank $r\leq n$
such that $H^{n-r+1}_\ast(\P^n,\kB)^\vee$ is a Cohen-Macaulay
$R$-module of dimension $\leq n-r$.\\[0.2cm]
This implies immediately:
\begin{corollary}
Let  $\kE$ be a (non-split) vector
bundle of odd rank $3\leq r \leq n$ on $\P^n$ with vanishing
intermediate cohomology,
with the exception of $H^{n-r+1}_{*}(\P^n,\kE)\neq 0$.
Let $s$ be a regular section of $\kE$. Then the top-dimensional
component of the $r$-codimensional zero-locus $Z(s)\subset\P^n$ is
arithmetically Gorenstein.
\end{corollary}

\begin{center} {\bf 4.~Relation to Good Determinantal Subschemes}
\end{center}
As an additional motivation for the study of Buchsbaum-Rim sheaves
$\kB_\phi$, we want to mention that the loci of regular sections 
of their duals $\kB^\ast_\phi$ actually correspond to a certain 
class of determinantal subschemes of projective space.
\begin{definition}
A subscheme $V=V_{t}(A)\subset {\P}^n$ is called determinantal 
if it is given as the zero set of all $t\times t$-minors of some 
homogeneous $g \times f$-matrix $A$  over $R=k[z_0,\dots,z_n]$. 
If $V_t$ has the expected codimension $(g-t+1)(f-t+1)$ it is 
called a standard determinantal subscheme.
\end{definition}
In codimension $2$ the situation is simple.
\begin{theorem}[Hilbert-Burch]
A codimension $2$ subscheme in ${\P}^n$ is standard determinantal if
and only if it is arithmetically Cohen-Macaulay.
\end{theorem}
Now let $A$ be a $t\times (t+r)$-matrix such that $V_{t}(A)$ has
the expected codimension $r+1$. Suppose that one can delete a 
generalized row from $A$ such that the ideal of maximal minors
of the resulting $(t-1)\times(t+r)$-matrix has the expected
codimension $r+2$. Then $V_{t}$ is called
a {\it good determinantal} subscheme.\\[0.2cm]
Every complete intersection in ${\P}^n$ is a good determinantal
subscheme. In parti\-cular, all codimension $2$ Gorenstein schemes
are good determinantal subschemes.\\[0.2cm]
{\bf Examples}\\[0.2cm]
The rational normal curve $V_{2}(A)\subset\P^n$ given by
$$ A=\left(
\begin{array}{ccccc}
z_0 & z_1 & z_2 & \dots & z_{n-1}\\
z_1 & z_2 & z_3 & \dots & z_n\\
\end{array}
\right)
$$
is a determinantal subscheme.
Let
$$ B=
\left(
\begin{array}{cccc}
z_1 & z_2 & z_3 & 0 \\
0  & z_1  & z_2 & z_3 \\
\end{array}
\right)
.$$
Then $V_{2}(B)\subset\P^3$ is standard determinantal but not good 
determinantal.\\[0.2cm]
The main result in the context with Buchsbaum-Rim sheaves is:
\begin{theorem}[\cite{Kreuzer}]
Let $X$ be a subscheme of ${\P}^n$ with $\mbox{codim}(X)\ge 2$. 
The following statements are equivalent :
\begin{itemize}
\item $X$ is a good determinantal subscheme of codimension $r+1$.
\item $X$ is a zero locus $Z(s)$ of a regular section 
      $s\in H^0(\P^n,\kB_\phi^\ast)$ of the {\bf dual} 
      $\kB^\ast_\phi$ of a Buchsbaum-Rim sheaf of rank $r+1$.
\item $X$ is standard determinantal and locally a complete
 intersection outside of some subscheme $Y\subset X$ of codimension
 $r+2$ in ${\P}^n$.
\end{itemize}
\end{theorem}
\begin{corollary}
A zero-scheme in ${\P}^3$ is good determinantal if and only 
if it is standard determinantal and a local complete intersection.
\end{corollary}
Let us mention that in the example above the ideal of $V_2(B)$ is
$(z_1,z_2,z_3)^2$. It is easy to see that it is not a local
complete intersection.

\begin{center} {\bf 5.~Implementation in Singular and Examples}
\end{center}
Now we want to describe how to put the construction method into practice.
The second author of this article implemented it in the computer
algebra system {\sl Singular}.\footnote{Singular is available at
{\tt http://www.singular.uni-kl.de/}} All the procedures in this and the
following section will soon be available as the Singular library 
{\tt buchsrim.lib}. Chris Peterson wrote similar procedures for 
Macaulay.\footnote{Macaulay 2 can be downloaded from
{\tt http://www.math.uiuc.edu/Macaulay2/}}\\[0.2cm]
The first question is: How can we compute a global
section of a Buchsbaum-Rim sheaf $\kB_{A}$ of rank $r$ in an 
algorithmical way?\\[0.2cm] 
We start with a homogeneous $t\times(t+r)$ matrix
$A$. Recall that $A$ is a block matrix with blocks $A_{ij}$ consisting of 
homogeneous polynomials in $R=k[z_0,\ldots,z_n]$. It can be considered 
as a part of the exact sequence
$$
0\lar \kB_{A} \lar \kF \stackrel{A}\lar \kG.
$$
Twist this sequence until we get global sections of $\kF$. Now apply
$H^{0}_{\ast}(\,\bullet\,)$. We get
$$
0 \lar H^{0}_{\ast}(\kB_{A}) \lar F \stackrel{A}\lar G.
$$
Let $H \stackrel{B}\lar F \stackrel{A}\lar G $ be a syzygy
sequence of $A$. The space of global sections of $\kB_{A}$ is just 
the degree zero component of $Im(B)$. So, we have to compute a matrix
representation $B$ of the first syzygy module of $A$, say via the 
command {\tt syz(A)}. Next, we 
take a linear combination $s=\sum_i\,f_i\,B_i$ of the columns $B_i$, 
where the coefficients $f_i$ are {\sl randomly chosen} homogeneous 
forms of some high degree $d$.
The ideal of the vanishing locus $Z(s)$ of the 
section is of course the transpose of the column vector $s$.\\
Given a homogeneous matrix $A$ and the degree $d$ of the forms $f_i$, 
the Singular procedure {\tt section(module A, int d)} generates such a
random section $s$ and returns the ideal of $Z(s)$.
One should check whether the section $s$ is regular, e.g.~by testing 
{\tt dim(std(section(A,d)))}$=n-r+1$.\\[0.2cm]
The most expensive part is the isolation of
the top-dimensional components $X$ of $Z(s)$. There are several ways 
to perform this computation: 
$$I_X=Ann_R Ext_R^r(\,R/I_{Z(s)},R\,)$$
This can be computed using the standard scripts
for $Ext$ and $Ann$ in Singular, but 
especially for high codimension $r$, this method is quite time-consuming:
One has to compute about $r$ Gr\"obner basis in order to
get $I_X$.\\ 
Instead, we use the following trick from liaison theory: Choose a 
regular sequence $J$ in $I_{Z(s)}$ of length $r$
(a randomly chosen sequence of homogeneous elements of high degree
is generically regular). Then the double ideal quotient
$(J:(J:I_{Z(s)}))$ is the saturated ideal $I_X$ 
(Singular procedure: {\tt top(ideal i)}).
This method is much faster, but one should definitely try to find 
another, more effective way to compute $I_X$.\\[0.2cm]
The procedure {\tt br($\ldots$)} includes everything described above
and returns the ideal of the top-dimensional part of a regular
section of the Buchsbaum-Rim sheaf $\kB_A$ associated to a homogeneous 
matrix $A$. This matrix can either be chosen randomly or be a specific one. 
As an option, one can obtain a detailed protocol of all the 
computations.\\[0.2cm]   
{\bf Example 2:} We give an annotated and slightly edited {\sl
  Singular} session which produces an arithmetically Gorenstein curve
of degree $21$ in $\P^6$.\\[0.2cm]
Consider a $1\times 6$-matrix $A$ of randomly chosen linear
forms on $\P^6$ such that 
\begin{equation} 0\lar\kB_A\lar
  6\,\kO_{\P^6}(2)\stackrel{A}{\lar}
\kO_{\P^6}(3)\lar \kM_A\lar 0
\end{equation}
is exact with $\mbox{Supp}(\kM_A)= Z$, $\dim Z=0$. The first
Chern class of
the reflexive rk-$5$ sheaf $\kB_A$ equals $9$. Using the Buchsbaum-Rim
resolution of $\kB_A$, we obtain 
$\deg S=\deg X = 21$.\\[0.2cm]
{\tt >ring r=32003,z(0..6),dp;\\
>LIB \"{}buchsrim.lib\"{};}\qquad\quad($\,$load the library 
{\tt buchsrim.lib}$\,$)\\
{\tt >ideal i=br(1,5,1,2);\\[0.2cm]
// We start with a random 1 x 6 matrix A of degree 1 forms:\\
//   A[1,1],A[1,2],A[1,3],A[1,4],A[1,5],A[1,6]\\[0.2cm]
// Check the codimension of Supp(Coker(A))...\\
// The codimension is 6 as expected.\\[0.2cm]
// The vanishing locus Z(j) of a randomly chosen section of the\\ 
// kernel sheaf is given by:\\[0.2cm]
// j[1]=z(0)*z(3)-12625*z(1)*z(3)+68*z(2)*z(3)+11333*z(3)\^{}2+} $\cdots$\\
{\tt // j[2]=z(2)\^{}2+15226*z(1)*z(3)+8747*z(2)*z(3)-12086*z(3)\^{}2+} 
  $\cdots$\\
$\vdots$\\
{\tt // j[6]=z(0)*z(1)-3269*z(1)*z(3)+9019*z(2)*z(3)+1243*z(3)\^{}2-}
$\cdots$\\[0.2cm]
{\tt // Its top-dimensional part Z(i) is given by:\\[0.2cm]
// i[1]=z(0)*z(3)-12625*z(1)*z(3)+68*z(2)*z(3)+11333*z(3)\^{}2+}
$\cdots$\\
{\tt // i[2]=z(2)\^{}2+15226*z(1)*z(3)+8747*z(2)*z(3)-12086*z(3)\^{}2+}
$\cdots$\\
{\tt // i[3]=z(1)*z(2)-7394*z(1)*z(3)+2944*z(2)*z(3)-4713*z(3)\^{}2-}
$\cdots$\\
{\tt // i[4]=z(0)*z(2)+309*z(1)*z(3)+2913*z(2)*z(3)+2392*z(3)\^{}2-}
$\cdots$\\
{\tt // i[5]=z(1)\^{}2-14366*z(1)*z(3)+11492*z(2)*z(3)+15380*z(3)\^{}2+} 
$\cdots$\\
{\tt // i[6]=z(0)*z(1)-3269*z(1)*z(3)+9019*z(2)*z(3)+1243*z(3)\^{}2-}
$\cdots$\\
{\tt // i[7]=z(0)*z(4)\^{}2-11212*z(1)*z(4)\^{}2+9631*z(2)*z(4)\^{}2-}
$\cdots$\\[0.2cm]
{\tt > hilb(std(i));}\quad ($\,$compute the Hilbert functions of a
standard basis of {\tt i}$\,$)\\[0.2cm]
{\tt //         1\quad\quad t\^{}0\\
//        -6 \quad t\^{}2\\
//        21 \quad t\^{}4\\
//       -21 \hspace{0.03cm}  t\^{}5\\
//         6 \hspace{0.37cm} t\^{}7\\
//        -1 \quad t\^{}9\\[0.2cm]
//         1 t\^{}0 } ($\,$ The coefficients of the $2^{\mbox{nd}}$
                            Hilbert function form the $h$-vector$\,$)\\  
{\tt//         5 t\^{}1\\
//         9 t\^{}2\\
//         5 t\^{}3\\
//         1 t\^{}4\\[0.2cm]
// codimension = 5} ($\,$Singular computes $\mbox{codim}_{\,\A^7} Z(I)$
   and  $\dim_{\,\A^7} Z(I)\,$)\\
{\tt // dimension   = 2\\
// degree      = 21} \qquad($\,\ldots$ as expected$\,$)\\[0.2cm]
From the $h$-vector $h_{I_X}=(1,5,9,5,1)$ we deduce 
the arithmetic genus $p_a(X)=20$, the Castelnuovo-Mumford 
regularity $\mbox{reg}(\kI_X)=5$ and $\kO_X\cong\omega_X(-2)$.
Using the command {\tt mres(i,0)}, we compute a {\sl minimal}
free resolution of the ideal $I_X$:
\begin{eqnarray*}
0\lar R(-9)\lar R(-6)\oplus\,6\,R(-7)\lar 21\,
   R(-5)\oplus\,R(-6)\lar & & \\ 
\lar R(-3)\oplus\,21\,R(-4)\lar 6\,R(-2)\oplus\,R(-3)\lar 
   I_X\lar 0 & &
\end{eqnarray*}
Convince yourself that this is exactly the resolution from 
Remark $3$ which is therefore minimal in this case.\\[0.2cm]
{\bf Example 3:} A Gorenstein threefold of degree $13$ in $\P^6$.\\[0.2cm]
Choose the matrix $A=(a_{11},a_{12},a_{13},a_{14})$ of quadratic forms
in the exact sequence 
$$0\lar \kB_A\lar 4\,\kO_{\P^6}(3)\stackrel{A}\lar\kO_{\P^6}(5)\lar\kM_A\lar 0$$
randomly and in such a way that the vanishing locus has the expected 
codimension $4$. Note that $c_1(\kB_A)=7$ and $\deg(X)=c_3(\kB_A)=13$. 
Again, we present an edited Singular session:\\[0.2cm]
{\tt >ring r=23,z(0..6),dp;\\
>LIB \"{}buchsrim.lib\"{};\\
>ideal i=br(1,3,2,3);}\\[0.2cm]
We compute the two Hilbert functions of a Gr\"obner basis of the ideal
{\tt i}$=I_X$.\\[0.2cm]
{\tt >hilb(std(i));}\\[0.2cm]
{\tt //         1 $\!\,$ t\^{}0\\
//        -1 t\^{}2\\
//        -4 t\^{}3\\
//         4 $\!\,$ t\^{}4\\
//         1 $\!\,$ t\^{}5\\
//        -1 t\^{}7\\[0.2cm]
//         1 t\^{}0\\
//         3 t\^{}1\\
//         5 t\^{}2\\
//         3 t\^{}3\\
//         1 t\^{}4\\[0.2cm]
// codimension = 3\\
// dimension   = 4}\qquad($\,$affine$\ldots\,$)\\
{\tt // degree      = 13}\\[0.2cm]
As expected, the arithmetically Gorenstein threefold $X$ has degree 
$13$. Its $h$-vector is $h_{I_X}=(1,3,5,3,1)$. Furthermore, we get
$\kO_X\cong\omega_X$ and $\mbox{reg}\,(\kI_X)=5$. 
A minimal free resolution of $I_X$ obtained via {\tt mres(i,0)} is
$$
0\lar R(-7)\lar 4\,R(-4)\,\oplus\,R(-5)\lar R(-2)\,\oplus\,4\,R(-3)\lar
I_X\lar 0,
$$
as stated in Theorem 2.
\begin{center} {\bf 6.~Application to Gorenstein Liaison} \end{center}
We shall apply the construction method also to deal with the 
following problem:\\
Let $V$ be equidimensional scheme in ${\P^n}$ of odd codimension
$c$.
We want to find an arithmetically Gorenstein subscheme $X$ of the same
codimension $c$, which contains $V$. Thus, we get a direct $G$-link 
$V\stackrel{X}{\sim} W$, where $W$ denotes the residue, i.e.~the
scheme associated to the saturated ideal $(\,I_X:I_V\,)$.
Surely, we can easily find a complete intersection with the desired
property. But from the point of view of Gorenstein liaison
we are interested in arithmetically Gorenstein schemes which are not 
complete intersections.\\[0.2cm]
{\bf Algorithm}
\begin{enumerate}
\item Choose a rk-$c$ Buchsbaum-Rim sheaf $\kB_{\phi}$ on ${\P^n}$.
$$
0\lar \kB_{\phi} \lar \kF \stackrel{\phi}\lar \kG \lar \kM_{\phi} \lar
0.
$$
\item Choose a regular section $s$ of $\kB_{\phi}(j)$ (for some shift $j
\in
{\mathbb Z}$), which is also in $H^{0}_{\ast}({\P^n}, \kF\otimes
\kI_{V})$.
Then the zero set $Z(s)$ of $s$ contains $V$.
\item Compute the top-dimensional part $X$ of the $Z(s)$. Then due to
Corollary 2, $X$ is arithmetically Gorenstein of codimension $c$ 
{\bf containing} $V$.
\end{enumerate}
{\bf Remark 5:}\\ 
We can always find a $j \in {\mathbb Z}$ such that
$H^0({\P^n}, \kB_{\phi}(j))\, \cap\, H^0({\P^n},\kF\otimes \kI_{V})
\ne 0.$\\[0.2cm]
Let us consider an example.\\[0.2cm]
{\bf Example 4:} Let $F$ be the Veronese surface in $\P^5$ given by the
$2\times 2$-minors of the matrix
\[
\left( \begin{array}{ccc} z_0 & z_3 & z_4\\
                           z_3 & z_1 & z_5\\
                           z_4 & z_5 & z_2
        \end{array} \right)
\]
Therefore $F$ is arithmetically Cohen-Macaulay (cf.~\cite{Curves},
p.~84). Consider the exact sequence
$$
0\lar \kB_{\phi} \lar 4\,\kO_{\P^5}(2)
\stackrel{\phi}{\lar} \kO_{\P^5}(3)
\lar \kM_\phi \lar 0,
$$
where $\phi:=(z_1,z_2,z_3,z_4)$ and the cokernel $\kM_\phi$ is
supported on the line
$L=Z(z_1,z_2,z_3,z_4)\subset\P^5$. The first syzygy module of $\phi$ is
$$
B=\left(
\begin{array}{cccccc}
0    &       0 &    0 &  -z_2 & -z_3   & -z_4\\
0    &    -z_3 & -z_4 &   z_1 &    0   &    0\\
-z_4 &     z_2 & 0    &    0  &    z_1 &    0\\
z_3  &    0    & z_2  &    0  &     0  &  z_1
\end{array}
\right).
$$
View the saturated ideal of the surface $F$
$$I_F=(\,z_3z_4-z_0z_5,\,z_1z_4-z_3z_5,\,z_2z_3-z_4z_5,\,z_1z_2-z_5^2,\,
  z_0z_2-z_4^2,\, z_0z_1-z_3^2\,)$$
as a $1 \times 6$-matrix. The global sections of
$4\,\kO_{\P^5}(2) \otimes \kI_{F}$ are linear combinations of columns of the
matrix
$I_F\otimes \mbox{Id}_{4}$. To find a common section $s$ of
$\kB_{\phi}(d)$ and
$4\,\kO_{\P^5}(2)\otimes \kI_{F}$, we consider the
intersection of the modules $B$ and $I_F\otimes\mbox{Id}_4$. Then $s$
is just a linear combination of the columns of the matrix
$B\cap( I_F\otimes\mbox{Id}_4)$ where
the coefficients are homogeneous $d$-forms in $z_0,\ldots, z_5$
whose degree $d$ has to be chosen large enough.\\
Computing this with Singular, one observes that the vanishing locus
$Z(s)$ of a random section $s$ is for example given by the ideal:
\begin{eqnarray*}
 I & = & ( z_1z_4-z_3z_5,\, z_2z_3-z_4z_5,\,z_0z_2-z_4^2,\, z_0z_1-z_3^2,\\
   &   &  z_3z_4^2-z_0z_4z_5,\, z_3^2z_4-z_0z_3z_5 )\;=\; I_F\cap
   (z_1,z_2,z_3,z_4)\cap (z_0,z_3,z_4).
\end{eqnarray*}
$Z(s)$ is consequently the union of $F$, containing the embedded line $L$,
and the plane $H=Z(z_0,z_3,z_4)$.
The regularity of the section $s$ can be easily checked by
computing the Hilbert function of a standard basis of $I$.
We get $\dim Z(s)=2$, $\deg Z(s)=5$ and
$h_I=(1,3,2,-1)$ as $h$-vector.\\
After isolating the top-dimensional part $X\subset Z(s)$,
we get a Gorenstein surface containing the Veronese
surface $F$. The saturated ideal of $X$ is
$$I_X=(\,z_3z_4-z_0z_5,\,z_1z_4-z_3z_5,\,z_2z_3-z_4z_5,\,
         z_0z_2-z_4^2,\,z_0z_1-z_3^2\,)=I_F\cap I_H.$$
The symmetric $h$-vector $h_{I_X}=(1,2,1)$ and the minimal free resolution
$$
0\lar R(-5)\lar 5\,R(-3)\lar 5\,R(-2)\lar I_X\lar 0
$$
confirm that $X$ is Gorenstein. The Veronese surface $F$ is therefore
$G$-linked to the plane $H$.   
\begin{center} {\bf 7.~Generalized Buchsbaum-Rim Sheaves} \end{center}
It is a natural question to ask what happens if the degeneracy locus
of $\phi\,:\,\kF\longrightarrow\kG$ does not have the expected
codimension $f-g+1$. Then the associated Buchsbaum-Rim complex is no
longer acyclic and the sections of the kernel sheaves often lose
the nice geometrical properties mentioned in Theorem $2$.\\[0.2cm] 
We restrict our attention to a particular class of examples where the 
degeneracy locus is an {\sl almost complete intersection}:\\
Let $G\subset\P^n$ be a codimension $3$ arithmetically Gorenstein
subscheme with $\kO_G\cong \omega_G(l)$. 
Choose a complete intersection $X$ of type $(d_1,d_2,d_3)$
containing $G$. Let $V$ denote the residue of $G$ under the CI-link,
i.e.~the variety associated to the ideal $\,(I_X:I_G)\,$, and let
$\alpha:=d_1+d_2+d_3$. Using the exact sequence
$$0\to\kI_X \to \kI_G\to\omega_V(l)\to 0,$$
the two resolutions 
$$
0\lar \kO_{\P^n}(-\alpha) \lar \bigoplus_{i<j}\kO_{\P^n}(-d_i-d_j)
 \lar \bigoplus_{i=1}^3 \kO_{\P^n}(-d_i)\lar \kI_X\lar 0
$$
and
$$
0\lar \kO_{\P^n}(l-n-1) \lar \kE_2 \lar \kE_1 \lar \kI_G \lar 0
$$
(\,$\kE_1$ and $\kE_2$ are decomposable bundles of the same rank $m$) 
and a mapping cone, we get
\begin{equation}\fbox{$0\lar\kK_{\phi} \lar\kF\stackrel{\phi}
{\lar}\kG\lar\kO_V(d)\lar
 0,$}\label{gebu}\end{equation}
where
$$\left.
\begin{array}{l}
\kF=\bigoplus_{i=1}^3\kO_{\P^n}(d-d_i)\oplus\kO_{\P^n}(n+1+d-l-\alpha)\\
\kG=\kO_{\P^n}(d)\\
\end{array}
\right\}\;d>\max_i\{d_i\}.$$
The degeneracy locus of those morphisms $\phi$ has therefore
codimension $3$ and not $4$ as expected.
Note that the bundles $\kF$ have global sections. 
\begin{definition} We call such kernel sheaves $\kK_\phi$ 
generalized Buchsbaum-Rim\\ sheaves.
\end{definition}
They have a free resolution
\begin{equation}0\lar
\kE_1^{\ast}(d-\alpha)\lar \bigoplus_{i<j}\kO_{\P^n}(d-d_i-d_j)\oplus\,
\kE_2^{\ast}(d-\alpha)\lar\kK_\phi\lar 0
\label{resgebu}\end{equation}
and are reflexive rk-$3$ sheaves (as 2nd syzygy sheaves).\\[0.2cm]
Our interest now focuses on properties of regular sections $s$ of generalized
Buchsbaum-Rim sheaves $\kK_\phi$. Unfortunately, their zero-loci $Z(s)$
are no longer arithmetically Gorenstein. However, one can still determine
a resolution of $I_{Z(s)}$. 
\begin{theorem}
Let $s$ be a regular section of a generalized Buchsbaum-Rim sheaf
$\kK_\phi$
on $\P^3$. Then the zero locus $Z(s)$ is an almost complete
intersection and
its saturated ideal $I_{Z(s)}\subset R=k[z_0,z_1,z_2,z_3]$ has a free 
resolution
\begin{eqnarray*}
0\lar R(\,\alpha-d-b\,)\oplus E^\ast_1(\,-b\,)\lar
\bigoplus_{i=1}^3 R(\,d_i-b\,)\oplus
\,R(\,-d\,)\oplus\,E_2^\ast(\,-b\,)\lar & &\\
\lar\bigoplus_{i=1}^3 R(\,d_i-d\,) \oplus\,R(\,d-b\,) \lar
I_{Z(s)}\lar 0, & &
\end{eqnarray*}
where $E_i:=H^0_\ast(\P^3,\kE_i)$ and $b:=2\,d-\alpha-l+4$.
\end{theorem}
{\bf Proof:} Let $\kA$ denote the image of $\phi$ in
(\ref{gebu}). Splitting the sequence into two short ones and applying
$\kH om(\,\bullet\,,\kO_{\P^3})$, we get  
$$0\lar\kH
om(\,\kO_V,\kO_{\P^3}(-d)\,)\lar\kG^{\ast}\lar\kA^{\ast}\lar
\kE xt^1(\,\kO_V,\kO_{\P^3}(-d)\,)\lar 0
$$
and
$$
0\lar\kA^{\ast}\lar\kF^{\ast}\lar\kK_{\phi}^{\ast}\lar
\kE xt^1(\,\kA,\kO_{\P^3}\,)\lar 0.
$$
$\kE xt^i_{\kO_{\P_3}}(\kO_V,\omega_{\P^3})$ vanishes for $0\leq
i<3$ and therefore $\kG^{\ast}\cong\kA^{\ast}$ implies
$\kA\cong\kO_{\P^3}(d)$. Hence the short exact sequence
$$0\lar\kK_{\phi}\lar\kF\lar\kA\lar 0$$
shows that $\kK_\phi$ is even locally free. Furthermore, we obtain 
a free resolution
$$0\lar G^\ast\lar F^\ast \lar H^0_{\ast}(\P^3,\kK_\phi^\ast)\lar 0$$
because $\kE xt^1(\,\kA,\kO_{\P^3}\,)=0$ and $H^1_{\ast}\kG=0$.
Let
\begin{equation}0\lar\kO_{\P^3}\stackrel{s}\lar\kK_{\phi}\lar\kC\lar 0
\label{secgebu}\end{equation}
be the exact sequence induced by the regular section $s$. Dualizing
it and using $\kE xt^1(\,\kK_\phi,\kO_{\P^3}\,)=0$, we get
$$0\lar\kC^\ast\lar\kK_{\phi}^{\ast}\lar\kI_{Z(s)}\lar 0.$$
Dualizing a second time, we realize that $\kC$ is a reflexive
rank $2$ sheaf. Therefore, $\kC^\ast\cong\kC(-c_1)$. 
Now use (\ref{resgebu}) and (\ref{secgebu}) in order to check
that $H^1_\ast(\,\P^3,\kK_\phi\,)\cong H^1_{\ast}(\,\P^3,\kC\,)=0$ and to
get a free resolution of $H^0_\ast\kC(-c_1)$ via a mapping cone:
\begin{eqnarray*}
0\rightarrow R(\,-c_1\,)\oplus E_1^\ast(\,d-\alpha-c_1\,) \rightarrow
\bigoplus_{i=1}^3 R(\,d+d_i-\alpha-c_1\,)\oplus
E^\ast_2(\,d-\alpha-c_1\,) & &\\ 
\rightarrow H^0_\ast\kC(-c_1)\rightarrow 0 & &
\end{eqnarray*}
Consequently, we obtain the following diagram:
\[\begin{array}{c@{\!\!}c@{\!\!}c@{\!}c@{\!\!}ccc}
    &         &      &     &  & 0 &  \\
    & R(-c_1) &      &\bigoplus_{i=1}^3 
    R(\,d+d_i-\alpha-c_1\,) & & \downarrow &  \\
  0\lar  & \oplus  & \lar & \oplus &
  \lar & H^0_\ast\kC(-c_1) & \lar 0\\
    &  E_1^\ast(\,d-\alpha-c_1\,) & &  E^\ast_2(\,d-\alpha-c_1\,) & & &  \\
    &         &      &  &
  &\big\downarrow &  \\
 & & & & & & \\
  0 \lar & G^\ast & \lar &   F^\ast & \lar & 
      H_\ast^0(\kK_\phi^\ast)&  \lar  0\\
 & & & & & \downarrow  &\\
 & & & & & I_{Z(s)} & \\
 & & & & & \downarrow & \\
 & & & & & 0. & 
\end{array}\]
Another application of the mapping cone lemma and the fact that
$c_1(\kC)=c_1(\kK_\phi)=\deg(\kF)-\deg(\kG)$ 
is equal to $3\,d-2\,\alpha-l+4=d-\alpha+b$ imply the claim. Note that $F$ is
a free $R$-module of rank $4$. Thus, $Z(s)$ is an almost complete
intersection of type $(d-d_1,d-d_2,d-d_3,d-\alpha-l+4)$.
\qed
{\bf Example 5:}
As Gorenstein scheme $G$, let us take $5$ points in $\P^3$ in general
position. Their saturated ideal $I_G$ can be obtained
using a regular section $s$ of the Buchsbaum-Rim sheaf
$$0\lar B_\phi \lar 4\,\kO_{\P^3}(2)\stackrel{\psi}\lar \kO_{\P^3}(3)
  \lar 0,$$
where $\psi=(\,z_0,z_1,z_2,z_3\,)$. Note that
$\kB_\psi=\Omega_{\P^3}^1(3)$ and $c_3(\kB_\psi)=5$. According to
Theorem $1$, there exists a free resolution
$$0\lar R\,(\,-5\,)\lar 5\,R\,(\,-3\,)\lar 5\,R\,(\,-2\,)\lar I_G\lar 0.$$ 
Now we choose a complete intersection $X$ of three cubic forms
containing $G$. The residue $V$ are $22$ points in $\P^3$. Their
saturated ideal $I_V$ is an almost complete intersection and has
the following minimal free resolution: 
$$
0 \lar 5\,R(\,-7\,)\lar 8\,R(\,-6\,)\lar 3\,R\,(\,-3\,)\oplus\,
R\,(\,-4\,)\stackrel{\phi} \lar I_V\lar 0.
$$
Compare it with (\ref{gebu}) and (\ref{resgebu}). Thus, the ``data'' for this
example is 
$$\kE_1=5\,\kO_{\P^3}(-2),\;\kE_2=5\,\kO_{\P^3}(-3),\;d_1=d_2=d_3=3,
\;\alpha=9\;\mbox{ and }\; l=-1.$$
We choose $d=6$. Using the Singular command {\tt
mres(section(syz($\phi$),3),0)} we get the a minimal free resolution
of the degree $13$ zero-locus $Z(s)$ of a regular section $s\in
H^0(\,\P^3,\kK_\phi\,)$:
$$
0\lar 4\,R(\,-6\,) \lar 7\,R(\,-5\,) \lar
3\,R(\,-3\,)\oplus\,R(\,-2\,)\lar I_{Z(s)}\lar 0$$
We recover exactly this sequence after deleting the ``ghost-summand''
$R(\,-5\,)\oplus\,R(\,-6\,)$ in the (non-minimal) free resolution in
Theorem $5$.\\[0.2cm]
{\bf Remark 6:} The authors believe that it is straightforward to
show that an analogon of Theorem $5$ holds 
for regular sections $s$ of generalized Buchsbaum-Rim sheaves
$\kK_\phi$ on $\P^n$, $n\geq 4$. One should use the fact that the 
Cohen-Macaulay type and the resolution are preserved under general 
hyperplane sections.

\vspace{1cm}
\begin{tabular}{lp{1.5cm}l}
Igor Burban & & Hans Georg Freiermuth\\
Universit\"at Kaiserslautern & & Columbia University \\
burban@mathematik.uni-kl.de & & freiermuth@math.columbia.edu
\end{tabular}
\end{document}